%% file: cap2.tex
\documentclass{article}

\usepackage{graphics}

\input mssymb12.tex
\input newmacro.tex

\title{Symplectic capacities of domains in $\Bbb C^2$}
\author{R. Hind \thanks{Supported in part by NSF grant DMS-0505778.}}

\date{\today}

\begin{document}

\maketitle

\section{Introduction}

In his paper \cite{gr} M. Gromov proved his celebrated
non-squeezing theorem. We will study domains $D$ in $\Bbb C^2$
with standard coordinates $(z_1,z_2)$ and projections $\pi_1$ and
$\pi_2$ onto the $z_1$ and $z_2$ planes respectively. The standard
symplectic form on $\Bbb C^2$ is $\omega = \frac{i}{2}
\sum_{j=1}^2 dz_j \wedge d\overline{z}_j$ and this restricts to a
symplectic form on the balls $B(r) = \{|z_1|^2 + |z_2|^2 < r^2
\}$. In this notation Gromov's non-squeezing theorem states that
if $\mathrm{area}(\pi_1(D)) \le C$ and there exists a symplectic
embedding $B(r) \to D$ then $\pi r^2 \le C$. Nowadays this can be
rephrased as saying that the Gromov width of $D$ is at most $C$.
Of course this is sharp when $D$ is a cylinder $\{|z_1| <r\}$.

For general $D$ it is natural to ask whether we can estimate the
Gromov width instead in terms of the cross-sectional areas
$\mathrm{area}(D \cap \{z_2=b\})$. But for any $\epsilon >0$ there
exists a construction of F. Schlenk, \cite{sch}, of a domain $D$
lying in a cylinder $\{|z_1|<1\}$ with Gromov width at least $\pi
- \epsilon$ but with all cross-sections having area less than
$\epsilon$. At least if we drop the condition on the domain lying
in the cylinder, the cross-sections can even be arranged to be
star-shaped, see \cite{book}. Nevertheless in this note we will
obtain such an estimate in terms of the areas of the
cross-sections for domains whose cross-sections are all starshaped
about the axis $\{z_1 =0\}$.

\begin{theorem}
Let $D \subset \Bbb C^2$ be a domain whose cross-sections $D \cap
\{z_2=b\}$ are star-shaped about center $z_1 =0$. Define $C=\sup_b
\mathrm{area}(\{z_2 =b\}\cap D)$. Then if $B(r) \to D$ is a
symplectic embedding we have $\pi r^2 \le C$. In other words, $D$
has Gromov width at most $C$.
\end{theorem}

In section $2$ we will establish an estimate on the Gromov width
for such domains $D$. This is combined with a symplectic embedding
construction to obtain our result in section $3$.

The author would like to thank Felix Schlenk for patiently
answering many questions.

\section{Embedding estimate}

Here we prove the following theorem.

\begin{theorem}
Fix constants $0 < K \le M$ and $0<t<1$. Let $D \subset \Bbb C^2$
be a domain of the form $D=\{r < c(\theta,z_2), |z_2|<M\}$ where
$(r,\theta)$ are polar coordinates in the $z_1$ plane and
$c(\theta,z_2)$ is a real-valued function satisfying $t \le
c(\theta,z_2) \le 1$ and $|\frac{\p c}{\p z_2}| \le \frac{1}{K}$.

Define $C=\sup_b \mathrm{area}(\{z_2 =b\}\cap D)$. Then if $B(r)
\to D$ is a symplectic embedding of the standard ball of radius
$r$ in $\Bbb C^2$ we have $\pi r^2 < C+3\sqrt{\frac{M}{tK^3}}$.
\end{theorem}

Its key implication for us is the following.

\begin{corollary}
Let $D=\{r < c(\theta,z_2), |z_2|<M\} \subset \Bbb C^2$ and
$C=\sup_c \mathrm{area}(\{z_2 =b\}\cap D)$. For any $L>0$ the
domain $D$ is a symplectic manifold with symplectic form $\omega_L
= \frac{i}{2}(dz_1 \wedge d\overline{z}_1+Ldz_2 \wedge
d\overline{z}_2)$. Let $r>0$ with $\pi r^2 >C$. Then for all $L$
sufficiently large the symplectic manifold $(D,\omega_L)$ does not
admit a symplectic embedding of the ball $B(r)$.
\end{corollary}

This follows by rescaling. Note above that the volume of
$(D,\omega_L)$ approaches infinity as $L \to \infty$.

{\bf Proof of Theorem $2$}

We consider the symplectic manifold $S^2 \times \Bbb C$ with a
standard product symplectic form $\omega = \omega_1 \oplus
\omega_2$ and still use coordinates $(z_1,z_2)$, where $z_1$ now
extends from $\Bbb C$ to give a coordinate on the $S^2 = \Bbb
CP^1$ factor. Still $\pi_1$ and $\pi_2$ denote the projections
onto the coordinate planes. Let $F$ be the area of the first
factor, we suppose that this is sufficiently large that the
complement of $\{z_1 =\infty\}$ can be identified with a
neighborhood of $\{|z_1| \le 1\}$ in $\Bbb C^2$, the
identification preserving the product complex and symplectic
structures. In other words, from now we assume that $D \subset S^2
\times \Bbb C \setminus \{z_1 = \infty\}$ and satisfies the
conditions on its cross-sections. Let $D^c$ denote the complement
of $D$ in $S^2 \times \Bbb C$.

Now let $\phi :B(r) \to D$ be a symplectic embedding. Then we
consider almost-complex structures $J$ on $S^2 \times \Bbb C$
which are tamed by $\omega$ and coincide with the standard product
structure on $D^c$. By now it is well-known, see \cite{gr}, that
for all such $J$ the almost-complex manifold $S^2 \times \Bbb C$
can be foliated by $J$-holomorphic spheres. In $\{|z_2| \ge M\}$
the foliation simply consists of the $S^2$ factors.

Let $S$ denote the image of the holomorphic curve in our foliation
passing through $\phi(0)$. By positivity of intersections $S$
intersects $\{z_1 = \infty\}$ in a single point, say $\{z_2=b\}$.
As above we will use polar coordinates $(r,\theta)$ in the plane
$\{z_2=b\}$. So we can write $D \cap \{z_2=b\} = \{r \le
c(\theta,b) := c(\theta)\}$. Let $A=\mathrm{area}(\{z_2 =b\}\cap
D)$. We intend to obtain lower bounds for both $\int_{S \cap D^c}
\omega_1$ and $\int_{S \cap D^c} \omega_2$.

First of all, we will suppose that $\pi_1(S \cap D^c)=\{r \ge
g(\theta)\}$ for a positive function $g$ and that $S \cap D^c$ is
a graph $\{z_2=u(z_1)\}$ over this region. We explain later how
essentially the same proof applies to the general case. Recall
that our assumptions imply that $t \le c(\theta),g(\theta) \le 1$
for all $\theta$. Define $h(\theta) = |g(\theta) - c(\theta)|$.

Define a holomorphic function $f:\{r \le \frac{1}{g(-\theta)}\}
\to \{|z_2| \le M\}$ by $f(z)=u(\frac{1}{z})$. Then $f(0)=b$ and
$|f(z)| \le M$ for all $z$. Therefore composing $f$ with a
translation we can redefine $f$ as a function $f:\{r \le
\frac{1}{g(-\theta)}\} \to \{|z_2| \le 2M\}$ with $f(0)=0$.

As $g(\theta) \le 1$ for all $\theta$ the map $f$ restricts to one
from $\{|z| \le 1\}$ and so by the Schwarz Lemma, if $|z|<1$ we
have $|f'(z)| \le \frac{2M}{1-|z|}$. On the boundary of the disk,
our assumptions on the boundary of $D$ imply that
$|f(\frac{1}{g(-\theta)}e^{i\theta})| \ge Kh(\theta)$.

Now we estimate \begin{eqnarray*} \int_{S \cap D^c} \omega_2 & = &
\mathrm{area}(\mathrm{image}(f))
\\
 & = & \int_0^{2\pi} d\theta \int_0^{\frac{1}{g(-\theta)}}
 r|f'(z)|^2 dr \\
 & = & \int_0^{2\pi} g(-\theta)d\theta
 \left( \int_0^{\frac{1}{g(-\theta)}}
 r|f'(z)|^2 dr \right) \left(\int_0^{\frac{1}{g(-\theta)}} dr \right) \\
 & \ge & t\int_0^{2\pi} d\theta \left(\int_0^{\frac{1}{g(-\theta)}
 }
 r^{\frac{1}{2}}|f'(z)| dr \right)^2.
 \end{eqnarray*}

 Now $$\int_0^{\frac{1}{g(-\theta)}
 } |f'(z)| dr \ge Kh(\theta)$$ and over all such functions
 $|f'(z)|$ the final integral above is minimized by taking
 $|f'(z)|$ as large as possible for small values of $r$. We
 compute $$\int_0^y \frac{2M}{1-r} dr =Kh(\theta)$$ when
 $y=1-e^{\frac{-Kh(\theta)}{2M}} < \frac{1}{g(-\theta)}$. Therefore putting $y=x^2$ we have

 \begin{eqnarray*}
t\int_0^{2\pi} d\theta \left(\int_0^{\frac{1}{g(-\theta)}}
 r^{\frac{1}{2}}|f'(z)| dr \right)^2 & \ge & t \int_0^{2\pi} d\theta \left( \int_0^{x^2}
 \frac{2M \sqrt{r}}{1-r} dr \right)^2 \\
  & = & 4M^2t \int_0^{2\pi} d\theta \left( \left[ -2\sqrt{r} +
  \ln\left( \frac{1+ \sqrt{r}}{1-\sqrt{r}}\right) \right]_0^{x^2}
  \right)^2 \\
  & = & 4M^2t \int_0^{2\pi} d\theta \left( -2x + \ln \left(
  \frac{1+x}{1-x} \right) \right)^2 \\
  & \ge & 4M^2t \int_0^{2\pi} \frac{4x^6}{9} d\theta
  \end{eqnarray*}

for the final estimate using the fact that $0<x<1$.

Now $$x^2 = 1-e^{\frac{-Kh(\theta)}{2M}} \ge
  (1-e^{-\frac{1}{2}})\frac{Kh(\theta)}{M}$$ since $\frac{Kh(\theta)}{2M}\le \frac{1}{2}$.

Therefore
  \begin{eqnarray*}
  \int_{S \cap D^c} \omega_2 & \ge & 4M^2t \int_0^{2\pi} \frac{4x^6}{9}
  d\theta\\
   & \ge & \frac{16}{9}(1-e^{-\frac{1}{2}})^3\frac{tK^3}{M} \int_0^{2\pi} h(\theta)^3 d\theta.
   \end{eqnarray*}

Next we compute
\begin{eqnarray*}
\int_{S \cap D^c} \omega_1 & = & F - \frac{1}{2} \int_0^{2\pi}
g(\theta)^2 d\theta \\
 & = & F -A - \frac{1}{2} \int_0^{2\pi}
(g(\theta)^2 - c(\theta)^2) d\theta \\
 & \ge & F -A - \frac{1}{2} \int_0^{2\pi}
(g(\theta) - c(\theta))(g(\theta) + c(\theta)) d\theta \\
 & \ge & F -A - \int_0^{2\pi} h(\theta) d\theta.
 \end{eqnarray*}

 Therefore writing $k=\frac{16}{9}(1-e^{-\frac{1}{2}})^3\frac{tK^3}{M}$ we
 have
 \begin{eqnarray*}
\int_{S \cap D^c} \omega & \ge & F-A - \int_0^{2\pi} (h(\theta) -
k h(\theta)^3) d\theta \\
 & \ge & F-A - 2\pi \frac{2}{3\sqrt{3k}} \\
 & = & F-A- \pi \sqrt{\frac{M}{3(1-e^{-\frac{1}{2}})^3tK^3}}.
 \end{eqnarray*}

 Thus $S \cap D$ has symplectic area at most $A+ \pi
 \sqrt{\frac{M}{3(1-e^{-\frac{1}{2}})^3tK^3}} < A + 3\sqrt{\frac{M}{tK^3}}$, since $S$ itself has area
 $F$.

We assumed above that $\pi_1(S \cap D^c)$ is starshaped about
$z_1=0$ and that $S \cap D^c$ is a graph over this region. If the
projection $\pi_1 :S \to \pi_1(S \cap D^c)$ is a branched cover
then we can define a function $f$ as before simply choosing a
suitable branch along the rays $\{\theta = \mathrm{constant}\}$.
The proof then applies as before. Now suppose that $\pi_1(S \cap
D^c)$ is not starshaped about $z_1=0$. Then we find the smallest
possible starshaped set $\{r \le g(\theta)\}$ containing the
complement of $\pi_1(S \cap D^c)$. The defining function $g$ will
then have discontinuities but this does not affect the proof which
again proceeds as before.

Finally we choose a $J$ which coincides with the push forward of
the standard complex structure on the ball $B(r)$ under $\phi$ but
remains standard outside $D$. The part of $S$ intersecting the
image of $\phi$ is now a minimal surface with respect to the
standard pushed forward metric on the ball and so must have area
at least $\pi r^2$, giving our inequality as required.

\section{Proof of Theorem $1$}

For any domain $E \subset \Bbb C^2$ we will write $C(E)=\sup_b
\mathrm{area}(\{z_2 =b\}\cap E)$. Again we let $C=C(D)$. Arguing
by contradiction suppose that $B(r) \to D$ is a symplectic
embedding with $\pi r^2
> C + \epsilon$.

Let $B$ be the image of the ball of radius $r$ in $D$. We will
prove Theorem $1$ by finding a symplectic embedding of $B$ into
$(D_1, \omega_L)$ for all sufficiently large $L$, where $D_1$ is a
domain $C^0$ close to $D$ and with $C(D_1) < C(D) + \epsilon$.
Such embeddings would contradict Corollary $3$.

First we choose a lattice of the $z_2$ plane sufficiently fine
that if we denote the gridsquares by $G_i$ then $\sup_i
\mathrm{area}(\pi_1(D \cap \pi_2^{-1}(G_i))) < C(D) + \epsilon$.
Then we let $D_1= \bigcup_i \pi_1(D \cap \pi_2^{-1}(G_i)) \times
G_i$, suitably smoothed.

Let $\{b_j\}$ be the vertices of our lattice. We make the
following simple observation.

\begin{lemma} Suppose that $B \cap \{z_2 = b_j\} = \emptyset$ for
all $j$. Then there exists a symplectic embedding of $B$ into
$(D_1, \omega_L)$ for all sufficiently large $L$.
\end{lemma}

{\bf Proof} It suffices to find a diffeomorphism $\psi$ of $\Bbb C
\setminus \{b_j\}$ which preserves the $G_i$ and such that
$\psi^*(L \omega_0) = \omega_0$, letting $\omega_0 = dz \wedge
d\overline{z}$ be the standard symplectic form. It is not hard to
construct such a map, and the product of this map on the $z_2$
plane with the identity map on the $z_1$ plane gives a suitable
embedding.

Given Lemma $4$, to find our embedding it remains to find a
symplectic isotopy of $D_1$ such that the image of $B$ is disjoint
from the planes $C_j = \{z_2 =b_j\}$. Equivalently we will find a
symplectic isotopy of the union of the $C_j$, compactly supported
in a neighborhood of $B$ and moving the $C_j$ away from $B$.

We may assume that the embedding of the ball of radius $r$ extends
to a symplectic embedding of a ball of radius $s$ where $s$ is
slightly greater than $r$. Let $U$ be the image of this ball and
$J_0$ the push-forward of the standard complex structure on $\Bbb
C^2$ to $U$ under the embedding.

\begin{lemma} There exists a $C^0$ small symplectic isotopy
supported near $\p U$ which moves each $C_j$ into a
$J_0$-holomorphic curve near $\p U$.
\end{lemma}

{\bf Proof} Let $(x+iy, u+iv)$ be local coordinates on $\Bbb C^2$.
Let $C$ be one of our curves. We may assume that in these
coordinates near to the origin $C \cap \p U$ is the curve
$\{(x,0,0,0)\}$ and therefore that nearby $C$ is the graph over
the $(x,y)$ plane of a function $h(x,y)=(u,v)$. So $u=v=0$ when
$y=0$.

There exists a constant $k$ such that $|u|$, $|v|$, $|\frac{\p
u}{\p x}|$ and $|\frac{\p v}{\p x}|$ are all bounded by $k|y|$
near $y=0$.

Now, such a graph is symplectic provided $$|\frac{\p u}{\p x}
\frac{\p v}{\p y} - \frac{\p v}{\p x} \frac{\p u}{\p y}| <1.$$

We can make $C$ holomorphic near $\p U$ by replacing $h$ by $(\chi
u, \chi v)$ where $\chi$ is a function of $y$, equal to $0$ near
$y=0$ and $1$ away from a small neighborhood. The resulting graph
remains symplectic provided $$|\chi \frac{\p u}{\p x}(\chi ' v+
\chi \frac{\p v}{\p y}) - \chi \frac{\p v}{\p x}( \chi 'u + \chi
\frac{\p u}{\p y})|<1$$ or rewriting
$$|\chi^2(\frac{\p u}{\p x}
\frac{\p v}{\p y} - \frac{\p v}{\p x} \frac{\p u}{\p y}) + \chi
\chi'(v\frac{\p u}{\p x} - u\frac{\p v}{\p x})|<1.$$

If we assume that $|\frac{\p u}{\p x} \frac{\p v}{\p y} - \frac{\p
v}{\p x} \frac{\p u}{\p y}|<1-\delta$ the graph remains symplectic
if $\chi$ is chosen such that $$|\chi \chi' (v\frac{\p u}{\p x} -
u\frac{\p v}{\p x})| < \delta$$ which is guaranteed if $\chi ' <
\frac{\delta}{ky^2}$.

Since the integral $\int_0^t \frac{\delta}{ky^2} dy$ diverges a
function $\chi$ satisfying this condition while being equal to $0$
near $0$ and $1$ away from an arbitrarily small neighborhood does
indeed exist as required. The resulting surface is clearly
isotopic through symplectic surfaces to the original $C$.

We now replace the $C_j$ by their images under the isotopy from
Lemma $5$. We let $J$ be an almost-complex structure on $U$ which
is tamed by $\omega$, coincides with $J_0$ near $\p U$, and such
that the $C_j \cap U$ are $J$-holomorphic.

Now $(U,J)$ is an (almost-complex) Stein manifold in the sense
that it admits a plurisubharmonic exhaustion function $\phi: U \to
[0,R)$. In fact, work of Eliashberg, see \cite{elf} and
\cite{mig}, implies that such a plurisubharmonic exhaustion exists
with a unique critical point, its minimum. Generically this will
be disjoint from the $C_j$.

Near the boundary we can take $\phi$ to be the push-forward under
the embedding of a function $\frac{|z|^N}{C}$ for some integer $N
\ge 2$ (depending perhaps on $U$) and (any given) constant $C$.
The definition of a plurisubharmonic function states that
$\omega_{\phi} = -dd^c \phi$ is a symplectic form on $U$ which is
compatible with $J$ (for a function $f$ we define $d^c f:= df
\circ J$). We can choose $C$ such that $\omega_{\phi} |_{\p U} =
\omega |_{\p U}$ and thus by Moser's lemma the symplectic
manifolds $(U, \omega)$ and $(U, \omega_{\phi})$ are
symplectomorphic via a symplectomorphism $F$ fixing the boundary.
In fact, adjusting the isotopy provided by Moser's method we may
assume that $F$ fixes the $C_j$ (since they are symplectic with
respect to both $\omega$ and $\omega_{\phi}$). Let $V$ denote the
image of $U \setminus B$ under $F$ and suppose that $\{\phi \ge
R_0\} \subset V$.

It now suffices to find a symplectic isotopy of the $C_j$ in
$(U,\omega_{\phi})$ moving the surfaces into the region $\{\phi
\ge R_0\}$. Then the preimages of these surfaces under $F$ gives a
symplectic isotopy moving them away from $B$ as required.

Let $Y$ be the gradient of $\phi$ with respect to the K\"{a}hler
metric associated to $\phi$. Equivalently $Y$ is defined by $Y
\rfloor \omega_{\phi} = -d^c \phi$. Define $\chi:[0,R) \to [0,1]$
to have compact support but satisfy $\chi(t)=1$ for $t \le R_0$.
Then the images of the $C_j$ under the one-parameter group of
diffeomorphisms generated by $X=\chi(\phi) Y$ will eventually lie
in $\{\phi \ge R_0\}$. Thus we can conclude after checking that
they remain symplectic during this isotopy. We recall that the
$C_j$ are $J$-holomorphic and finish with the following lemma.

\begin{lemma} Let $G$ be a diffeomorphism of $U$ generated by the
flow of the vectorfield $X$. Then $G^* \omega_{\phi}(Z,JZ)>0$ for
all non-zero vectors $Z$.
\end{lemma}

{\bf Proof} For any function $f$ we compute
\begin{eqnarray*}
{\cal L}_X f(\phi)d^c \phi & = & f'(\phi) X \rfloor d\phi \wedge
d^c \phi + f(\phi)X \rfloor dd^c \phi + d(f(\phi) X \rfloor d^c
\phi) \\
 & = & (f'(\phi)d\phi(X) + f(\phi)\chi(\phi))d^c \phi.
 \end{eqnarray*}

Thus $G^* d^c \phi = g(\phi)d^c \phi$ for some function $g$ and
$$G^* \omega_{\phi} = g(\phi)\omega_{\phi} - g'(\phi)d\phi \wedge
d^c \phi.$$

The function $g$ is certainly positive and so $G^* \omega_{\phi}$
evaluates positively on the (contact) planes $\{d \phi = d^c \phi
=0\}$. Therefore if $G^* \omega_{\phi}$ evaluates nonpositively on
a $J$-holomorphic plane then there exists such a plane containing
$Y$. But this is clearly not the case, as $G^*
\omega_{\phi}(Y,JY)=\omega_{\phi}(G_*Y,G_*JY)=-k d^c \phi(G_* JY)$
for some positive constant $k$ and $-d^c \phi(G_* JY)=-G^* d^c
\phi (JY) = g(\phi) d \phi(Y)>0$.

\pagebreak

Richard Hind\\
Department of Mathematics\\
University of Notre Dame\\
Notre Dame, IN 46556\\
email: hind.1@nd.edu

\end{document}

%% file: newmacro.tex
\def\NABLA#1{{\mathop{\nabla\kern-.5ex\lower1ex\hbox{$#1$}}}}
\def\Nabla#1{\nabla\kern-.5ex{}_#1}

\newcommand{\p}{{\partial}}

\newcommand{\qed}{\hfill$\Box$\medskip}
\newtheorem{theorem}{Theorem} 
\newtheorem{corollary}[theorem]{Corollary}
\newtheorem{lemma}[theorem]{Lemma}

%% file: cap2.bbl
\begin{thebibliography} {99}

\bibitem{elf} Y. Eliashberg, Filling by holomorphic disks and its applications, {\it London Math. Society Lecture Notes}, Series 151 (1991), 45-67.

\bibitem{mig} Y. Eliashberg, Symplectic Geometry of plurisubharmonic functions, with notes by Miguel Abreu, NATO Adv. Sci. Inst. Ser. C Math. Phys. Sci., 488, {\it Gauge theory and symplectic geometry (Montreal, PQ 1995), 49-67, Kluwer Acad. Publ., Dordrecht, 1997}.

\bibitem{gr} M. Gromov, Pseudo-holomorphic curves in symplectic manifolds, {\it Inv. Math.}, 82(1985), 307-347.
 \bibitem{sch} F. Schlenk, On a question of Dusa McDuff. {\it Int. Math. Res. Not.}, 2003, no. 2,
 77--107.
 \bibitem{book} F. Schlenk, Embedding problems in symplectic geometry, de Gruyter Expositions in Mathematics, 40, Berlin, 2005.


\end{thebibliography}
